\documentclass[11pt,twoside]{book}
\usepackage[centertags]{amsmath}
\usepackage{amsfonts}
\usepackage{amssymb}
\usepackage{amsthm}
\usepackage{newlfont}
\usepackage{color}
\usepackage{graphicx}

\date{\color{black}  2013 November}
\def \n {\noindent}
\usepackage{fancyhdr}
\usepackage[final]{pdfpages}
\usepackage[T1]{fontenc} 
\usepackage{textcomp}
\usepackage{marvosym}
\usepackage{pdfpages}

\begin{document}

\begin{center}
{\large {\bf {\color{red}On a generalization of the Perron-Frobenius theorem \\in an ordered Banach space}}}\\
\end{center}

\begin{center}
{\bf Abdelkader Intissar}$^{1,2}$
\end{center}

\n $^{1}$ Equipe d'Analyse Spectrale, Université de Corse, UMR-CNRS No. 6134, Quartier Grossetti\\
20 250 Corté, France\\
Tél : 00 33 (0) 4 95 45 00 33 -Fax : 00 33 (0) 4 95 45 00 33\\
Email address : intissar@univ-corse.fr\\

\n $^{2}$ Le Prador, 129 rue du Commandant Rolland, 13008 Marseille, France\\
Email address : abdelkader.intissar@orange.fr\\

\n {\bf {\color{red} Abstract}}\\

\n We work in an ordered Banach space with closed generating positive cone. We show that a positive compact operator has zero spectral radius or a positive eigenvector with the  corresponding eigenvalue equal to the spectral radius.\\

\n {\bf {\color{red} $\S$ 1  The Perron-Frobenius Theorem in $\mathbb{R}^{n}$}}\\

\n {\color{blue} {\bf  (A) Nonnegative Vectors and Matrices}}\\

\n {\bf Definition 1.1}\\

\n {\color{red}$\bullet$} A vector $x  \in \mathbb{R}^{n}$ is nonnegative, and we write $x  \geq  0$, if its coordinates
are nonnegative. It is positive, and we write $x > 0$, if its coordinates are (strictly) positive. Furthermore, a matrix $A \in  M_{n×m}(\mathbb{R})$ (not necessarily square) is nonnegative (respectively, positive) if its entries are nonnegative (respectively, positive); we again write $A \geq  0$ (respectively, $A > 0$). More generally, we define an order relation $x \leq  y$ whose meaning is $y- x \geq 0$.\\

\n {\color{red}$\bullet$} Given $x  \in \mathbb{C}^{n}$, we let $|x|$ denote the nonnegative vector whose coordinates are the numbers $\mid x_{j}\mid$ Likewise, if $A \in M_{n}(\mathbb{C})$, the matrix $|A|$ has entries $|a_{i j}|$.\\

\n{\color{red}$\bullet$} Observe that given a matrix and a vector (or two matrices), the triangle inequality implies $|Ax| \leq |A| · |x|$.\\

\n For a systematic study of positive matrices,  we can consult the volumes 1 and 2 of Gantmacher {\bf{\color{blue}[Gantmacher]}}.\\

\n {\bf Proposition 1.2}\\

\n  A matrix is nonnegative if and only if $x \geq 0$ implies $Ax \geq 0$. It is positive if and only if $x \geq 0$ and $x \neq 0$ imply $Ax > 0.$\\

\n {\bf Proof}\\

\n  Let us assume that $Ax \geq 0$ (respectively,$ > 0$) for every $x \geq 0$ (respectively, $\geq 0$ and  $\neq 0$). Then the $ith$ column $A^{(i)}$ is nonnegative (respectively, positive), since it is the image of the $ith$ vector of the canonical basis. Hence $A \geq  0$ (respectively, $> 0$).\\

\n Conversely, $A \geq 0$ and $x \geq 0$ imply trivially $Ax \geq 0$. If $A > 0, x \geq  0$, and $x \neq 0$, there exists an index $l$ such that $x_{l} > 0$. Then
$\displaystyle{(Ax)_{i} = \sum_{j}a_{ij}x_{j } \geq  a_{il}x_{l} > 0}$, and hence $Ax > 0$. \\

\n {\bf Primitive and irreducible non-negative square matrices}\\

\n {\bf Definition 1.3} (see  for example {\bf{\color{blue}[Sternberg]}})\\

\n {\color{red}$\bullet$} A non-negative matrix square $A$ is called {\color{red}primitive} if there is a $k$ such that all the entries of $A^{k}$ are positive. \\

\n  {\color{red}$\bullet$} It is called {\color{red}irreducible} if for any $i, j$ there is a $ k = k(i, j) $ such that $\displaystyle{(A^{k})_{ij} > 0}$.\\

 \n {\color{red}$\bullet$} An $n \times n $ matrix $A = (a_{ij}) $ is said to be reducible if $n \geq 2$ and there exists a permutation matrix P such that :\\
\n $^{t}PAP = \left [ \begin{array} {cc} A_{11}&A_{12}\\
\quad\\
0&A_{22}\\
\end{array} \right ]$
\quad\\

\n where  $A_{11}$ and $A_{12}$ are square matrices of order at least one. If A is not reducible, then it is said to be irreducible.\\

\n {\bf Lemma 1.4}\\

\n If $A$ is irreducible then $I + A$ is primitive.\\

\n {\bf Proof} \\

\n  Indeed, the binomial expansion $\displaystyle{I + A)^{k} =  I + kA + \frac{k(k-1)}{2}A^{2} + ...}$ will eventually have positive entries in all positions if $k$ large
enough.\\

\n An important point is the following:\\

\n {\bf Proposition 1.5}\\

\n If $A \in M_{n}(\mathbb{R}) $ is nonnegative and {\color{red}irreducible}, then $\displaystyle{(I+A)^{n-1} > 0}$.\\

\n {\bf Proof}\\

\n  Let $x \neq 0$ be nonnegative, and define $x^{m} = (I +A)^{m}x$, which is nonnegative too. Let us denote by $P_{m}$ the set of indices of the nonzero components of $x^{m}$:
$P_{0}$ is nonempty. Because $\displaystyle{x_{}^{m+1} \geq x_{i}^{m}}$ , one has $\displaystyle{P_{m} \subset P_{m+1}}$. Let us assume that the cardinality $|P_{m}|$ of $P_{m}$ is strictly less than $n$. There are thus one or more zero components, whose indices form a nonempty subset $I$, complement of $P_{m}$. Because $A$ is {\color{red}irreducible}, there exists some nonzero entry $a_{ij}$, with $i \in I $and $j \in P_{m}$. Then $\displaystyle{x_{i}^{m+1} \geq a_{ij}x_{i}^{m} > 0}$, which shows that $P_{m+1}$ is not equal to $P_{m}$, and thus $\displaystyle{|Pm+1| > |Pm|}$.\\

\n By induction, we deduce that  $\displaystyle{|P_{m}| \geq  min{m+1,n}}$. Hence $\displaystyle{|P_{n-1}| = n}$, meaning that $\displaystyle{x_{n-1} > 0}$. We conclude with Proposition 1.2.\\

\n  {\color{blue} {\bf  (B) The Perron-Frobenius Theorem:Weak Form}}\\

\n Here we denote by $\sigma(A)$ the spectrum (the set of all eigenvalues) of a (square) matrix $A$, and by $\rho(A)$ the spectral radius of $A$, i.e., the quantity $max \{ \mid \lambda \mid ;  \lambda \in \sigma(A)\}$\\

\n {\bf Theorem 1.6} \\

\n Let $A \in M_{n}(\mathbb{R})$ be a nonnegative matrix. Then its spectral radius $\rho(A)$ is an eigenvalue of $A$ associated with a nonnegative eigenvector.\\

\n {\bf Proof}\\

\n  Let $\lambda$ be an eigenvalue of maximal modulus and $v$ an eigenvector, normalized by $\mid\mid v \mid\mid_{1} = 1$. Then $\displaystyle{\rho(A) \mid v|\mid = |\lambda v| = |Av| \leq  A|v|}$.\\

\n Let us denote by $\mathcal{C}$ the subset of $\mathbb{R}^{n}$ (actually a subset of the unit simplex $\mathcal{K}_{n}$) defined by the (in)equalities $\displaystyle{ \sum_{i}^{} x_{i} = 1, x \geq 0}$, and $Ax \geq \rho(A)x$. This is a {\color{red}closed convex set}, nonempty, inasmuch as it contains $|v|$. Finally, it is bounded, because $x \in \mathcal{C}$ implies $0 \leq x_{j} \leq 1$ for every $ j$; thus it is {\color{red} compact}. \\

\n Let us distinguish {\color{red}two cases}.\\

\n {\color{blue}1} There exists $x \in  \mathcal{C}$ such that  $Ax = 0$. Then $\rho(A)x \leq 0$ furnishes {\color{red}$\rho(A) = 0$}. The theorem is thus proved in this case.\\

\n {\color{blue}2} For every $x \in \mathcal{C}, Ax  = 0$. Then let us define on $\mathcal{C} $ a {\color{red}continuous map} $f$ by\\

\n $\displaystyle{f (x) = \frac{1}{\mid\mid Ax \mid\mid_{1}} Ax}$.\\

\n It is clear that $f (x) \geq 0$ and that $\mid\mid f (x) \mid\mid_{1} = 1$. Finally,\\

\n $\displaystyle{Af (x) = \frac{1}{\mid\mid Ax \mid\mid_{1}} A Ax = \frac{1}{\mid\mid Ax \mid\mid_{1}} A\rho(A)x = \rho(A)f(x) }$.\\

\n so that $ f (\mathcal{C}) \subset \mathcal{C}$. Then {\color{red}Brouwer's theorem} (see {\bf{\color{blue}[Berger et al]}}, p. 217) asserts that a continuous function from a {\color{red}compact convex subset }of $\mathbb{R}^{n}$ into itself has a fixed point.\\

\n Thus let $y$. be a fixed point of $f$. It is a nonnegative eigenvector, associated with the eigenvalue $r = \mid\mid Ay \mid\mid_{1}$. Because $y \in \mathcal{C}$, we have $ry = Ay \geq \rho(A)y$ and thus $r \geq \rho(A)$, which implies $r = \rho(A)$.\\

\n {\bf Remark 1.7}\\

\n (i) That proof can be adapted to the case where a real number $r$ and a nonzero vector $y$ are given satisfying $y \geq 0$ and $Ay \geq ry$. \\

\n Just take for $\mathcal{C}$ the set of vectors $x$ such that $\displaystyle{\sum_{i}^{ }x_{i} = 1, x \geq  0}$, and $Ax \geq rx$. We then conclude that $\rho(A) \geq r$.\\

\n (ii) In section 2 (the main of this work), we prove that a positive compact operator has either zero spectral or a positive eigenvector with the corresponding eigenvalue equal to the spectral radius.\\

\n  {\color{blue} {\bf (C) The Perron-Frobenius Theorem: Strong Form}}\\

\n {\bf Theorem 1.8}\\

\n  Let $A \in M_{n(}\mathbb{R})$ be a {\color{red}nonnegative irreducible }matrix. Then $\rho(A)$ is a {\color{red}simple } eigenvalue of $A$, associated with a {\color{red}positive}eigenvector. Moreover, $\rho(A) > 0$.\\

\n {\bf Proof}\\

\n For $r \geq 0$, we denote by $\mathcal{C}_{r}$ the set of vectors of $\mathbb{R}^{n}$ defined by the conditions \\

\n $x \geq 0$, \quad $\mid\mid x \mid\mid_{1} = 1$ and  $Ax \geq rx$.\\

\n Each $\mathcal{C}_{r }$ is a {\color{red}convex compact} set. We know  that if $\lambda$ is an eigenvalue associated with an eigenvector $x$ of unit norm $\mid\mid x \mid\mid_{1} = 1$, then $\displaystyle{\mid x \mid \in \mathcal{C}_{\mid \lambda\mid}}$. In particular, $\mathcal{C}_{\rho(A)}$ is nonempty.\\

\n Conversely, if $\mathcal{C}_{r }$ is nonempty, then $ x \in \mathcal{C}_{r}$, $r = r\mid\mid x \mid\mid_{1} \leq \mid\mid Ax \mid\mid _{1}  \leq A \mid\mid  x \mid\mid _{1} = \mid\mid A \mid\mid _{1}$, and therfore $ r \leq \mid\mid A \mid\mid _{1}$.\\

\n  Furthermore, the map $r  \longrightarrow \mathcal{C}_{r }$ is non increasing with respect to inclusion, and is `` left continuous'' in the following sense. If $ r > 0$, one has $\displaystyle{\mathcal{C}_{r} = \bigcap_{s < r} \mathcal{C}_{r}}$.\\

\n Let us then define $\displaystyle{R = sup\{r ; \mathcal{C}_{r}  \neq \emptyset\}}$, so that $R \in  [\rho(A),  \mid\mid A \mid\mid_{1}]$. The monotonicity with respect to inclusion shows that $r < R$ implies $\mathcal{C}_{r}  \neq \emptyset$.\\
\n If $x > 0$ and $\mid\mid x \mid\mid_{1} = 1$, then $Ax > 0 $ because $A$ is {\color{red} nonnegative} and {\color{red} irreducible}.\\

\n Setting $\displaystyle{ r := min_{j}(Ax)_{j}/x_{j} > 0}$, we have $\mathcal{C}_{r}  \neq \emptyset$, whence $R \geq r > 0$. The set $\mathcal{C}_{R}$, being the intersection of a totally ordered family of nonempty compacts sets, is nonempty.\\

\n Let $x \in\mathcal{C}_{R}$ be given. the lemma 1.9 below \\

\n {\bf Lemma 1.9}\\

\n Let $r \geq 0$ and $x \geq 0$ such that $Ax \geq rx$ and $Ax \neq rx$. Then there exists $r' > r$ such that $\mathcal{C}_{r'}$ is nonempty.\\

\n shows that $x$ is an eigenvector of $A$ associated with the eigenvalue $R$.We observe that this eigenvalue is not less than $\rho(A)$ and infer that $\rho(A) = R$. Hence $\rho(A)$ is an eigenvalue associated with the eigenvector $x$.\\

\n The following lemma :\\

 \n {\bf Lemma 1.10}\\

\n The nonnegative eigenvectors of $A$ are positive. The corresponding eigenvalue is positive too.\\

\n ensures that $x > 0$ and $\rho(A) > 0$.\\

 \n The simplicity of the eigenvalue $\rho(A)$ is given in\\
 
 \n {\bf Lemma 1.11}\\
 
 \n The eigenvalue $\rho(A)$ is simple.\\
 
 \n Finally, we can state the following result:\\
 
 \n {\bf Lemma 1.12}\\
 
\n  Let $M, B \in M_{n}(\mathbb{C}) $ be matrices, with $M$ irreducible and $|B| \leq M$. Then \\

\n $\rho(B) \leq \rho(M)$.\\

\n In the case of equality ($\rho(B) = \rho(M)$), the following hold:\\

\n {\color{red}$\bullet$} $|B| = M$.\\

\n {\color{red}$\bullet$} For every eigenvector $x$ of $B$ associated with an eigenvalue of modulus $\rho(M)$, $|x|$ is an eigenvector of $M$ associated with $\rho(M)$.\\

\n {\color{blue}{\bf (D) Proof of Lemmas}}\\

\n {\bf Proof of Lemma 1.9}\\

\n Set $y := (I_{n} + A)^{n-1}x$. Because $A$ is irreducible and $x \geq 0$ is nonzero, one has $y > 0$. Likewise, $Ay - ry = (I_{n} + A)^{n-1}(Ax - rx) > 0$. 

\n Let us define $r := min_{j}(Ay)_{ j}/y_{ j}$, which is strictly larger than $r$. We then have $Ay \geq ry$, so that $\mathcal{C}_{r}$ contains the vector $y/\mid\mid y \mid\mid_{1}$.\\

\n {\bf Proof of Lemma 1.10}\\

\n Given such a vector $x$ with $Ax = \lambda x$, we observe that $\lambda \in \mathbb{R}_{+}$. Then $\displaystyle{x = \frac{1}{(1+ \lambda)^{n-1}}(I_{n} + A)^{n-1}x}$ and  the right-hand side is strictly positive, from Proposition 1.5. Inasmuch as $A$ is irreducible and nonnegative, we infer $Ax  = 0$. Thus $\lambda = 0$; that is, $\lambda > 0$.\\

\n {\bf Proof of Lemma 1.11}\\

\n Let $P_{A}(X)$ be the characteristic polynomial of $A$. It is given as the composition of an $n$-linear form (the determinant) with polynomial vector-valued functions
(the columns of $XI_{n} - A)$. If $\phi$ is $p$-linear and if $\displaystyle{V_{1}(X), . . . ,V_{p}(X) }$ are polynomial vector-valued functions, then the derivative of the polynomial $\displaystyle{P(X) := \phi(V_{1}(X), . . . ,V_{p}(X)) }$ is given by \\

\n $\displaystyle{P'(X) = \phi(V_{1}^{'}, V_{2}, . . . ,V_{p})+ \phi(V_{1}, V_{2}^{'}, . . . ,V_{p}) +· · ·+\phi(V_{1}, V_{2}, . . . ,V_{p}^{'})}$.\\

\n One therefore has \\

\n $\displaystyle{P_{A}(X) = det(e_{1},a_{2}, . . . ,a_{n}) + det(a_{1},e_{2}, . . . ,a_{n})+· · ·+det(e_{1},a_{2}, . . . ,a_{n-1}, e_{n})}$,\\

\n where $a_{j}$ is the $jth$ column of $XI_{n} - A$ and $ \{e_{1}, . . . ,e_{n} \}$ is the canonical basis of $\mathbb{R}^{n}$.\\

\n Developing the $jth$ determinant with respect to the $jth$ column, one obtains\\

\n $\displaystyle{P_{A}^{'}(X) =\sum_{j=1}^{n}P_{A_{j}}(X)}$\\

\n where $A_{j} \in M_{n-1}(\mathbb{R})$ is obtained from $A$ by deleting the $jth$ row and the $jth$ column.\\

\n Let us now denote by $B_{j} \in M_{n}(\mathbb{R})$ the matrix obtained from $A$ by replacing the entries of the $jth$ row and column by zeroes. This matrix is block-diagonal, the two diagonal blocks being $A_{j } \in M_{n-1}(\mathbb{R})$ and $0 \in M_{1}(\mathbb{R})$. Hence, the eigenvalues of $B_{j}$ are those of $A_{j}$, together with zero, and therefore $\rho(B_{j}) = \rho(A_{j})$. Furthermore, $|B_{j}| \leq A$, but $|B_{j}|  = A$ because $A$  is irreducible and $B_{j}$ is block-diagonal, hence reducible. It follows (Lemma 12) that $\rho(Bj) < \rho(A).$ Hence $P_{A_{j }}$ does not vanish over $[\rho(A),  +\infty)$. Because $P_{A_{j}} (t) \sim t^{n-1}$ at infinity, we deduce that $P_{A_{j }}(\rho(A)) > 0.$\\

\n Finally, $P_{A}(\rho(A))$ is positive and $\rho(A)$  is a simple root.\\

\n {\bf Proof of Lemma 1.12}\\

\n In order to establish the inequality, we proceed as above. If $\lambda$ is an eigenvalue of $B$, of modulus $\rho(B)$, and if $x$ is a normalized eigenvector, then $\rho(B)|x| \leq |B| · |x| \leq M|x|$, so that $\mathcal{C}_{\rho(B)}$ is nonempty. Hence $\rho(B) \leq R = \rho(M)$\\

\n Let us investigate the case of equality. If $\rho(B) = \rho(M)$, then $ |x| \in \mathcal{C}_{\rho(M)}$, and therefore $|x|$ is an eigenvector: $M|x| = \rho(M)|x| = \rho(B)|x| \leq |B| · |x|$. Hence, $ (M- |B|)|x| \leq 0$. Because $ |x| > 0$ (from Lemma 10) and $M-|B| \geq 0$, this gives $|B| = M$.\\

\n The above results  have a long history, in fact  In 1907, Perron  {\bf {\color{blue}[Perron1]}} and  {\bf {\color{blue}[Perron2]}} gave proofs of the following famous theorem, which now bears his name, on positive matrices:\\

\n {\color{black}{\bf Perron Theorem (1907)}}\\

\n Let $A$ be a square positive matrix. Then $\rho(A)$ is a simple eigenvalue of $A$ and there is a corresponding positive eigenvector.\\
Furthermore,    $\mid \lambda \mid < \rho(A)$ for all $\lambda \in \sigma(A), \lambda \neq \rho(A)$.\\

\n \n In $1912$, Frobenius  {\bf {\color{blue}[Frobenius]}} extended the Perron theorem to the class of irreducible nonnegative matrices:\\

\n {\color{black}{\bf Frobenius Theorem (1912)}}\\

\n Let $A \geq 0$ be irreducible. Then \\

\n  (i)  $\rho(A)$ is simple eigenvalue of $A$, and there is a corresponding positive eigenvector.\\

\n  (ii) If $A$ has $m$ eigenvalues of modulus $\rho(A)$, then they are in the following form $\displaystyle{\rho(A)e^{\frac{2ik\pi}{m}}}$ ; \\
 \quad \quad \quad  $ k = 0, ...., m-1$.\\

\n   (iii) The spectrum of $A$ is invariant under a rotation about the origin of the complex plane\\
   \quad \quad \,\quad  by $\frac{2\pi}{m}$, i.e., $\displaystyle{e^{\frac{2i\pi}{m}}\sigma(A) = \sigma(A)}$.\\

\n  (iv) If $m > 1$ then there exists a permutation matrix $P$ such that :\\
\n $^{t}PAP = \left (\begin{array}{ccccc} 0&A_{12}&&&\\
\quad\\
&0&A_{23}&&\\
\quad\\
 & &0&\ddots&\\
\quad\\
&&&\ddots&A_{m-1,m}\\
\quad\\
A_{m,1}&&&&0\\
\end{array} \right )$
\quad\\
\n where the zero blocks along the diagonal are square.\\

\n we refer the reader to the intersting paper of C. Bidard and M. Zerner  {\bf {\color{blue}[Bidard et al]}} for an application of the Perron-Frobenus theorem in relative spectral theory and {\bf{\color{blue}[Alintissar et al]}} on some economic models.\\

\n A natural extension of the concept of a nonnegative matrix is that of an integral operator with a nonnegative kernel. The following extension of Perron's theorem is due to Jentzsch {\bf{\color{blue}[Jentzsch]}}:\\

\n {\color{black}{\bf Jentzsch Theorem (1912)}}\\

\n Let $k(., .)$ be a continuous real function on the unit square with $k(s, t) > 0$ for all $0 \leq s, t \leq 1$. If $K: L^{2}[0, 1] \longrightarrow  L^{2}[0, 1]$ denotes the integral operator with kernel $k$ defined by setting \\

\n $\displaystyle{ (Kf)(s) = \int_{0}^{1}k(s, t)f(t)dt , f \in L^{2}[0, 1],}$\\

\n then\\

\n (i) $K$ has positive spectral radius;\\

\n (ii) the spectral radius $\rho(K)$ is a simple eigenvalue, with (strictly) positive eigenvector;\\

\n (iii) if $ \lambda = \rho(K)$ is any other eigenvalue of $K$, then $\mid \lambda \mid < \rho(K)$.\\

\n A generalization of Jentzsch theorem is given by Schaefer in his book {\bf {\color{blue} [Schaefer]}} as follow:\\

\n {\color{black}{\bf Schaefer Theorem(1974)}}\\

\n  Let $(\mathbb{X}, \mathcal{T}, \tau)$ be a measure space with positive measure $\tau$ and $L_{p}(\tau) $ be the set of all measurable functions on $\mathbb{X}$ whose absolute value raised to the $p^{-th}$ power has finite integral.\\

\n Let $T$ be an integral bounded operator defined on $L_{p}(\tau) $ by a kernel $\mathcal{N} \geq 0$.\\

\n We suppose that:\\

\n (i) There exists $n \in \mathbb{N}$ such that $^{n}$ is compact operator.\\

\n (ii) For $\mathbb{S} \in \mathcal{T}; \tau(\mathbb{S} > 0$ and $\tau(\mathbb{X} - \mathbb{S}) > 0$ we have :\\

\n $\displaystyle{\int_{\mathbb{X} - \mathbb{S}}\int_{\mathbb{S}}\mathcal{N}(s, t)d\tau(s)d\tau(t) > 0}$ \\

\n Then the spectral radius $r(T)$ of integral operator $T$ is simple eigenvalue associated to an eigenfunction $f$ satisfying $f(s) > 0 \, \tau$-almost everywhere.\\

\n An fine study  on above theorems  is given in an article of Zerner in 1987 {\bf {\color{blue}[Zerner]}} in particular the following result : \\

 \n {\color{black}{\bf Zerner Theorem (1987)}}\\

\n Suppose that $A$ is irreducible and that $\rho(A)$ is a pole of its resolvent. Then $\rho(A)$ is non-zero and it is a simple pole, any positive eigenvector associated with $\rho(A)$ is quasi-interior and any positive eigenvector of the transpose $A'$ of $A$ associated with $\rho(A)$ is a form strictly positive.\\

\n  If moreover $\rho(A)$ is of finite multiplicity, it is a simple eigenvalue. \\

\n {\bf Remark 1.13}\\

\n {\color{red}$\bullet$} We should not be under any illusions about the conclusion `` $ \rho (A) $ nonzero ''. In applications, we do not really see how we could prove that $ \rho (A) $ is a pole sarees at the same time that it is non-zero. \\

\n \n {\color{red}$\bullet$} We refer the reader to an original application of Jentzsch Theorem in reggeons field theory  see T. Ando and M. Zerner $(1984)$ in  {\bf {\color{blue}[Ando et al]}} and A. Intissar  and J.K Intissar $(2019)$  in {\bf {\color{blue}[Intissar et al]}} \\

\n {\bf {\color{red} $\S$ 2  The Perron-Frobenius Theorem in an ordered Banach space with closed generating positive cone}}\\

  \n {\bf {\color{blue} (A) Introduction and statement of the results}}\\
  
\n  Below we give some  abstract properties of  ordered relation on a Banach space\\

\n Let $(E, \leq )$ be an ordered vector space and  $\mathcal{C}$  is a proper $\mathcal{C}\cap(-\mathcal{C}) = \{0\}$ convex positive cone in $E$,\\

\n  {\color{blue}$\bullet$}  $\displaystyle{x \leq y \Longrightarrow  x +a \leq y+a }$\\

\n  {\color{blue}$\bullet$} $\displaystyle{x \leq y \iff y - x \in \mathcal{C}}$\\

\n {\color{blue}$\bullet$} $\mathcal{C} + \mathcal{C} \subset \mathcal{C}$\\

\n  {\color{blue}$\bullet$}  $\displaystyle{x \leq y \, and \, t \geq 0 \Longrightarrow  tx  \leq ty }$\\

\n  {\color{blue}$\bullet$}  $\displaystyle{x_{1} \leq y_{1} \, and \, x_{2} \leq y_{2}    \Longrightarrow x_{1} + y_{1} \leq x_{2} + y_{2}}$\\

\n  {\color{blue}$\bullet$}  $\displaystyle{ x \in \mathcal{C} \, and \,  -x \in \mathcal{C} \Longrightarrow x = 0}$\\

\n {\color{blue}$\bullet$}  $\displaystyle{\mathcal{C}}$ , is  closed.\\
 
\n {\bf Definition 2.1}\\

\n A real Banach space $E$ equipped with a closed convex proper and generating cone $\mathcal{C}$ is called an ordered Banach space.\\

\n By proper is meant that there exists no $x \neq 0$ in $E$ such that both $x$ and $-x$ belong to $\mathcal{C}$ and by generating that  every $x \in E$ has a decomposition  $x = x^{+} - x^{-}$ where $x^{+} $ and $x^{-}$ both belong  to $\mathcal{C}$\\

\n The cone $\mathcal{C}$ defines an order relation in $E$.\\

\n A vector in $E$ is called positive (agreement with the French terminology) if it belong to $\mathcal{C}$ , negative if its opposite is positive.\\

\n The elements belonging to the cone $\mathcal{C}$ are called positive elements of $E$\\

\n An operator  $A$ acting on an ordered Banach space $E$ is said to be positive if it transforms positive elements into positive elements i.e. $A(\mathcal{C}) \subset \mathcal{C}$.\\

 \n In $1948$, in an abstract order-theoretic setting, in the important memoir {\bf{\color{blue}[14]}} Krein and Rutman have partially extended the Perron-Frobenius theorem to a  positive compact linear operator leaving invariant a convex cone in a Banach space.\\

\n They obtained the following:\\
 
\n {\color{black}{\bf Krein-Rutman Theorem (1948)}} \\

\n (i) Let $A$ be a positive compact linear operator on  $E$. Suppose that $A(\mathcal{C}) \subseteq \mathcal{C}$, where $\mathcal{C}$ is a closed generating cone in
$E$. If $\rho(A) > 0$, then there exists a nonzero vector $x \in \mathcal{C}$ such that $Ax = \rho(A)x$.\\

\n (ii) Let $A$ be a positive compact linear operator on  $E$. Assume there is a non negative vector $x$, a natural number $p$ and a positive real $\lambda$ such that \\

\n $\displaystyle{A^{p}x \geq \lambda^{p}x}$ \hfill { } {\bf {\color{blue} (2.1)}}\\

\n (The inequality $x \leq y$ means that $y - x \in \mathcal{C}$)\\

\n Then $A$ has a positive eigenvector associated with an eigenvalue at least equal to $\lambda$.\\

\n Remember that the spectral radius of an operator is the radius of the smallest closed disk containing its spectrum. In the case of a compact operator, it is either zero or the largest modulus of its eigenvalues.\\

\n This suggest the following result which is the main part of this work:\\

\n {\bf Theorem 2.2}\\

\n The spectral radius $\rho(A)$ of a positive compact linear operator $A$ on $E$ is either $0$ or an eigenvalue corresponding to a positive eigenvector.\\

\n {\bf Proof}\\

\n We give the proof in several steps.\\

\n This section will be devoted to proof of this theorem. But before, we give the following corollary and remark and two classical theorems from Kato's book on the structure of spectrum of {\color{red}compact } operator\\

\n {\bf Corollary 2.3} (monotonicity of the spectral radius)\\

\n Let $A_{1}$ and $A_{2}$ be  two positive compact linear operators on $E$ such that $A_{2} - A_{1}$ is positive.\\

\n Then the spectral radius of $A_{2}$ is at least equal to the spectral radius of $A_{1}$.\\

\n {\bf Proof}\\

\n Call $r$ the spectral radius of  $A_{1}$.  We may assume  $r {\color{red} >} 0$, otherwise the conclusion  is obvious. By theorem 2.2, $r$ is an eigenvalue of $A_{1}$ corresponding  to a positive eigenvector $u$. Then we have $A_{2}u \geq A_{1}u = r u$.\\

\n So that by (ii) of Perron Frobenius theorem with $p = 1$ $A_{2}$ has an eigenvalue at least equal to $r$.\\

\n {\bf Remark 2.4}\\

\n From now on, $A$ will denote a positive compact linear operator on $E$ with {\color{red}non zero spectral radius}.\\

\n Multiplying $A$ by a positive factor multiplies its spectral radius by the same factor. So we can assume the spectral radius of $A$ to be equal to one and we shall do so henceforth.\\

\n A simple and well know special case is if we suppose that $A$ have an eigenvalue of the form $\displaystyle{e^{\frac{ik\pi}{n}}}$. Then one is an eigenvalue of $A^{2n}$.\\

\n Let $u$ be an eigenvector associated with this eigenvalue $(A^{2n}u = u)$. \\

\n We assume $u$ to be non negative otherwise we would take $-u$ instead. So we can apply (ii) of Krein-Rutman theorem, obtaining a positive eigenvector of $A$ associated with an eigenvalue at least equal to {\color{red}one} but  it cannot be larger.\\

\n The spectrum of a {\color{red}compact} operator $A$ in $E$ has a simple structure ``analogous'' to that of an operator in a finite-dimensional space. This is translated by the following theorem\\

\n {\bf Theorem 2.5} (Kato theorem III.6.26 p. 185, {\bf {\color{blue}[Kato]}})\\

\n Let $A$ be a compact operator acting on $E$. Then \\

\n (i)  its spectrum $\sigma(A)$ is a cuntable set with no accumulation point different from zero.\\

\n (ii) each nonzero $\lambda \in \sigma(A)$ is an eigenvalue of $A$ with finite multiplicity,  and $\overline{\lambda}$ is an eigenvalue of $A'$ with same multiplicity.\\

\n In this section we will also use the following theorem adapted from theorem III. 6.17 p. 178 [Kato] on the separation of spectrum.\\

\n {\bf Theorem 2.6} (Kato theorem III.6.17 p. 178, {\bf {\color{blue}[Kato]}})\\

\n $ E$ can to be split into the direct sum of two closed subspaces $E'$ and $E''$, both invariant under $A$, with the following properties where $A'$ and $A''$ denote the restriction of $A$ to $E'$ and $E''$ respectively.\\

\n (i) All eigenvalues of $A'$ have modulus one.\\

\n (ii) The spectral radius of $A''$ is at most equal to $r$ (equal if we have chosen $r$ as small as possible).

\n Moreover, the eigenvalues of modulus one being isolated and of finite multiplicity, $E'$ has a finite dimension.\\

\n It follows from above theorem that a crucial step in the proof of theorem 2.2 is that $E'$ contains a {\color{red} non zero positive vector}.\\

\n  This will be proved  after the following  proposition in following step:\\
 
  \n {\bf {\color{blue} (B) Construction of a sequence of nearly eigenvectors}}\\
  
  \n {\bf Proposition 2.7}\\
  
  \n The aim of this proposition is to construct a sequence of vectors $w_{k}$ and a sequence $p_{k}$ with the following properties:\\
  
  \n ($\alpha$)  $w_{k}$ is positive and has norm one.\\
  
  \n ($\beta$) $\displaystyle{lim \, p_{k} = + \infty}$ as $k \longrightarrow \infty$.\\
  
  \n ($\gamma$) There is  a sequence of real numbers $\lambda_{k}$ such that :\\
  
  \n (i) $\displaystyle{lim \, inf \lambda_{k}  \geq 1}$ as $k \longrightarrow \infty$.\\
  
  \n (ii)  $\displaystyle{lim \, z_{k} = 0}$ as $k \longrightarrow \infty$. by setting $\displaystyle{z_{k} = \lambda_{k}w_{k} - A^{p_{k}}w_{k}}$.\\
  
  \n {\bf Proof}\\
  
   \n As a preliminary, we need the following simple consequence of Hahn-Banach . A positive form on $E$ is a form which is non negative on $\mathcal{C}$ :\\
  
  \n {\bf Lemma 2.8}\\
  
  \n Let $u$ be any non negative vector.  then there is a continuous positive linear form $f$ on $E$ such that {\color{red}$f(u) > 0$}.\\
  
  \n {\bf Proof of lemma}\\
  
  \n $- \mathcal{C}$ is a closed convex set and $u \notin (- \mathcal{C})$. So there is an {\color{red} affine} form  $f_{1}$ which is {\color{red}non positive} on ({\color{red}$- \mathcal{C}$}) and {\color{red} positive} at $u$ ({\color{blue}{\bf [Bourbaki] }}ch. 11 $\S$ 3 proposition 4)\\

  \n Let us define $f$ by :\\
  
  \n $\displaystyle{ f_{1}(x) = f(x) + f_{1}(0)}$. \hfill { } {\bf {\color{blue} (2.2)}}\\
  
  \n and check that it has the desired properties.\\
  
  \n First $\displaystyle{ f(u) = f_{1}(u) - f_{1}(0)}$ where $f_{1}(u)$ is positive and $f_{1}(0)$ is not, thus  $f(u) > 0$.\\
  
  \n Now let $x \in \mathcal{C}$ and look at   $\displaystyle{ f_{1}(tx) = tf(x) + f_{1}(0)}$ which we know to be non positive for all negative $t$ implying $f(x) \geq 0$. The lemma is proved.\\
  
  \n Now for technicalities, we may assume that $A$ has an eigenvalue $e^{i\theta}$ where $\frac{\theta}{\pi}$ is irrational. Let $u + iv$. be a corresponding eigenvector. Here again we may and will assume $u$ non negative and we write  $v = v^{+} - v^{-}$ with  $v^{+}$ and $ v^{-}$ positive.\\
  Let $f$ given  by the above lemma, we define an operator $B$ by \\
  
  \n $\displaystyle{B(x) = \frac{f(x)}{f(u)}v^{+}}$ \hfill { } {\bf {\color{blue} (2.3)}}\\
  
  \n {\bf Remark 2.9}\\
  
  \n Note that $B$ is positive, compact (even of rank one) and maps $u$ on $v^{+}$.\\

 \n  Now  let $p$ be any natural number such that $p \in ]0, \frac{\pi}{2}[ \, mod \, 2\pi $ then we have :\\
  
  \n $\displaystyle{A^{p}(u + iv) = e^{ip\theta}(u + iv)}$. \hfill { } {\bf {\color{blue} (2.4)}}\\
  
  \n Taking the real pats of above equality to get \\
  
  \n $\displaystyle{A^{p}u = cos (p\theta)u - sin (p\theta) v}$ \hfill { } {\bf {\color{blue} (2.5)}}\\
  
  \n So that\\
  
  \n $\displaystyle{(A^{p} +  sin (p\theta)B)u = cos (p\theta)u + sin (p\theta) v^{-1} \geq cos (p\theta) u}$ \hfill { } {\bf {\color{blue} (2.6)}}\\
  
  \n We apply the (ii) of Krein-Rutman theorem to the operator $\displaystyle{A^{p} +  sin (p\theta)B}$, the number $p$ of the statement of the theorem being here equal to one.\\
  
 \n  We conclude that there is a positive vector $x_{p}$ of norm one such that:\\
 
\n  $\displaystyle{(A^{p} +  sin (p\theta)B)x_{p} = \mu_{p}x_{p}}$ \hfill { } {\bf {\color{blue} (2.7)}}\\
 
 \n with\\
 
 \n $\displaystyle{ \mu_{p} \geq cos(p\theta)}$ \hfill { } {\bf {\color{blue} (2.8)}}\\
 
 \n As $\displaystyle{\frac{\theta}{\pi}}$ is irrational, we can find a strictly increasing sequence $(p_{k})$ with :\\
 
 \n $\displaystyle{p_{k}\theta = \epsilon_{k} \quad mod. \,  [2\pi]}$, $\epsilon_{k} > 0$ , $\displaystyle{lim \, \epsilon_{k} = 0} $ as $k \longrightarrow \infty$  \hfill { } {\bf {\color{blue} (2.9)}}\\
 
 \n Setting $\displaystyle{w_{k} = x_{p_{k}}}$ and $\displaystyle{\lambda_{k} = \mu_{p_{k}}}$ we check ($\alpha$) and ($\beta$) and we have:\\
 
 \n $\displaystyle{z_{k} = sin (p_{k}\theta)Bw_{k}}$ \hfill { } {\bf {\color{blue} (2.10)}}\\
 
 \n so that\\
 
 \n $\displaystyle{\mid\mid z_{k} \mid\mid \leq \epsilon_{k}\mid\mid B \mid\mid}$ \hfill { } {\bf {\color{blue} (2.11)}}\\
 
\n  and by (2.8):\\
 
 \n $\displaystyle{\lambda_{k} \geq cos(p_{k}\theta)}$ \hfill { } {\bf {\color{blue} (2.12)}}\\
 
 \n where the second member has limit one when $k$ tends to infinity. Then properties ($\alpha$) to  ($\gamma$) hold.\\
 
 \n In this step, we consider the projecteur $P'$ with kernel $E''$ and range $E'$ and let $P'' = I - P'$ then we will show the following lemma :\\
 
 \n {\bf Lemma 2.10}\\
 
\n $\displaystyle{ lim \, P''w_{k} = 0 }$ as $k \longrightarrow \infty$\\

\n {\bf Proof} \\

\n By the property ($\gamma$) of preceding step, we have $\displaystyle{\lambda_{k}w_{k} = A^{p_{k}}w_{k} + z_{k}}$ where $\lambda_{k}$ is bounded from below and $z_{k}$ converges to zero. So what we have to show is that $\displaystyle{P'' A^{p_{k}}w_{k}}$ converges  to zero.\\

\n Let $r'$ be a number satisfying  $r < r' <1$  (As the non zero eigenvalues of $A$ are isolated , there is a number  $r < 1$ such that the modulus of any eigenvalue of $A$ is either one or not larger than $r$).\\

\n Notice that \\

\n $\displaystyle{ P'' A^{p_{k}}w_{k} =  A^{p_{k}}P''w_{k} = A''^{p_{k}}P''w_{k}}$  \hfill { } {\bf {\color{blue} (2.13)}}\\

\n By Gel-fand's theorem, we know that the spectral radius of $A''$ is the limit of $\displaystyle{ \mid\mid A''^{n} \mid\mid^{\frac{1}{n}}}$ soi that we have, for $k$ large enough, the following inequality :\\

\n $\displaystyle{\mid\mid A''^{p_{k}}\mid\mid < r'^{p_{k}}}$  \hfill { } {\bf {\color{blue} (2.14)}}\\

\n Now, putting together  (2.13) and  (2.14) we get: \\

\n $\displaystyle{ \mid\mid P'' A^{p_{k}}w_{k}\mid\mid \leq  r'^{p_{k}} \mid\mid P'' \mid\mid}$  \hfill { } {\bf {\color{blue} (2.15)}}\\

\n As $r'$ is smaller than one, this ends the proof of the lemma.\\

\n {\bf Remark 2.11} (Extracting a convergent subsequence) \\

\n $\displaystyle{(P'w_{k})}$ is a bounded sequence in the finite dimensional vector space $E'$. We can extract from it a convergent subsequence and, by above lemma,  the corresponding subsequence of the $w_{k}'s$ converges to the same limit $w$.\\

\n  As a limit of the $\displaystyle{P'w_{k}'s}$ , $w$ belongs to $E'$. As a limit of the $w_{k}'s$, it is positive and has norm one.\\

\n In this last step we will see the situation  in $E'$ in particular we will derive the existence of a positive eigenvector in $E'$ by the following Lemma:\\

 \n {\bf Lemma 2.11}\\  
 
 \n Let $\mathcal{C}$ be a closed convex cone in finite dimensional vector space $F$. Assume   $\mathcal{C}$ to be  neither $\{0\}$ nor the whole space.\\
 
 \n Let $T$ be a linear operator on $F$ which maps  $\mathcal{C}$ into itself. \\
 
 \n Then $T$ has an eigenvector belong to  $\mathcal{C}$.\\
 
 \n {\bf Proof} \\
 
 \n If $T$ maps some non zero vector of $\mathcal{C}$ to $0$, we are through. If not, call $\mathbb{S}$ the unit sphere of some euclidian metric on $F$ and $\mathcal{C}_{1}$ the intersection of $\mathbb{S}$ and $\mathcal{C}$.\\
 
 \n $\mathcal{C}_{1}$ is homeomorphic to closed ball of some space $\mathbb{R}^{k}$ (take for instance the stereographic projection from some point of $\mathbb{S}$ not belong to  $\mathcal{C}$ ; we get a compact convex set). \\
 
 \n The mapping $\displaystyle{x \longrightarrow \frac{T(x)}{\mid\mid T(x) \mid\mid}}$ is continuous from $\mathcal{C}_{1}$ into itself. By Brouwer's therem {\bf {\color{blue}[Brouwer]}}, it has a fixed point, and this is an eigenvector of $T$. The lemma is proved.\\
 
 \n {\bf Summing up}\\
 
 \n We have proved the existence of a positive eigenvector in $E'$. As has already been indicated in step 2, the corresponding eigenvalue is positive, the eigenvector being positive; it has modulus one, the eigenvector belong to $E'$. It must be one, the spectral radius of $A$.\\
 
 \n {\bf References}\\
 
 \n {\bf {\color{blue} [Alintissar et al]}}  Alintissar, A. , Intissar, A. and Intissar, J.K. :. On dynamics of wage-price spiral and stagflation in some
model economic systems,  adsabs.harvard.edu/abs/2018arXiv181201707A\\
 
  \n {\bf{\color{blue}[Ando et al] }} Ando, T. and Zerner, M. :. sur une valeur propre d'un opérateur, Commun. Math. Phys. 93,123 139 (1984)\\
 
\n {\bf{\color{blue} [Berger et al]}}, Berger , M. : and Gostiaux, B. :. .Differential Geometry: Manifold, Curves and Surfaces, volume 115 of Graduate Texts in Mathematics. Springer-Verlag, New York, 1988.\\
 
\n {\bf{\color{blue} [Bidard et al]}} Bidard,C and Zerner M. :. The Perron-Frobenius theorem in relative spectral theory ?, Mathematische Annalen, vol.289, (1991) pp.451-464.\\

\n {\bf{\color{blue}[Bourbaki] }} Bourbaki, N. :. Eléments de mathématique livre V , Espaces vectoriels topologiques Hermann (1953)\\

\n {\bf{\color{blue}[Frobenius]}}  Frobenius, G.F. :. Uber Matrizen aus nicht negativen Elementen, Sitzungsber. Kon. Preuss.
Akad. Wiss. Berlin, (1912), 456-477\\
 
 \n {\bf{\color{blue}[Gantmacher1]}} Gantmacher., F. R. :. The Theory of Matrices. Vol. 1. Chelsea, New York, 1959.\\
 
 \n {\bf{\color{blue}[Gantmacher2]}} Gantmacher, F. R. :. The Theory of Matrices. Vol. 2. Chelsea, New York, 1959.\\
 
  \n {\bf{\color{blue}[Kato]}} Kato, T. Perturbation Theory for linear Operators Springer second Edition (1980)\\
 
\n {\bf{\color{blue} [Krein-Rutman]}} Krein, M.G and Rutman , M. A. :. Linear operators leaving invariant a cone in a Banach space, Amer. Math. Soc. Transl. Ser. 1 10 (1950), 199-325 [originally Uspekhi Mat. Nauk 3 (1948), 3-95].\\

\n {\bf{\color{blue} [Intissar et al]}} Intissar, A. and  Intissar, J.K. :.  A Complete Spectral Analysis of Generalized Gribov-Intissar's Operator in Bargmann Space, Complex Analysis and Operator Theory (2019) 13:1481?1510 \\

\n {\bf{\color{blue} [Jentzsch]}}  Jentzsch, P. :. Uber Integralgleichungen mit positivem Kern, J. Reine Angw. Moth. 141 (1912). 235-244.\\

\n {\bf{\color{blue}[Perron1]}} Perron, O. :. Grundlagen fur eine Theorie des Jacobischen Kettenbruchalgorithmus, Math. Ann. 63 (1907), 1-76.\\

\n {\bf{\color{blue}[Perron2]}} Perron, O. :. Zur Theorie der uber Matrizen, Math. Ann. 64 (1907), 248-263.\\

\n {\bf{\color{blue}[Schaefer]}}  Scheafer, H. H. :.Banach Lattices and Positive Operators, Springer, Berlin/Heidelberg/ New York, (1974).\\

\n {\bf{\color{blue}[Sternberg]}} Sternberg, S. :. The Perron-Frobenius theorem, Lecture 12,\\
{\color{blue}http://www.math.harvard.edu/library/sternberg/slides/1180912pf.pdf}\\

\n {\bf{\color{blue}[zerner] }} Zerner, M. Quelques proprités spectrales des opérateurs positifs, journal of Functional Analysis, 72, 381-417 (1987)\\

  \end{document}